\documentclass[runningheads,a4paper]{llncs}
\usepackage[utf8]{inputenc}
\usepackage{hyperref}
\usepackage{apacite}
\usepackage{verbatim}
\usepackage{mathtools}
\usepackage{color}
\usepackage{amssymb,amstext,amsmath}
\newcommand{\overbar}[1]{\mkern 1.5mu\overline{\mkern-1.5mu#1\mkern-1.5mu}\mkern 1.5mu}
\newcommand{\no}[1]{\overbar{#1}}

\def\F{\mathcal F}

\def\pr{\mathbb{P}}
\def\prev{\mathbb{P}}

\def\S{\mathcal{S}}
\def\C{\mathcal{C}}

\renewcommand{\bot}{\emptyset}
\renewcommand{\top}{\Omega}

\def\ind{\vbox{\hbox{$\bot$\kern-.6em$\bot$}\kern-.05em}}
\def\dep{\vbox{\hbox{$\top$\kern-.6em$\top$}\kern-.05em}}
\newcommand*{\normally}{\mathrel{\ooalign{$|$\hfil\cr\kern+1pt$\thicksim$}}} 
\newcommand*{\nnormally}{\mathrel{\ooalign{$|$\hfil\cr\kern+1pt$\thicksim$}\negthickspace \negthickspace /} } 

\def\S{\mathcal{S}}
\def\C{\mathcal{C}}

\def\F{\mathcal{F}}

\title{Probabilistic entailment and iterated conditionals}
\author{
Angelo Gilio\inst{1} \thanks{Retired}   
	\and Niki Pfeifer\thanks{Supported by his DFG project PF~740/2-2 (within the SPP1516)} \inst{2} 
	\and Giuseppe Sanfilippo\inst{3} 
}
\institute{
		Department SBAI,
	University of Rome ``La Sapienza'', Italy
	\email{angelo.gilio@sbai.uniroma1.it}  
	\and
	Munich Center for Mathematical Philosophy, LMU Munich, Germany
 \email{niki.pfeifer@lmu.de} 
 \and
 	Department of Mathematics and Computer Science,
 University of Palermo, Italy  
 \email{giuseppe.sanfilippo@unipa.it}
}
\renewcommand\top{\Omega}
\renewcommand\bot{\emptyset}
\begin{document}
	\maketitle
\begin{abstract}
In this paper we exploit the notions of conjoined and iterated conditionals, which are defined in the setting of coherence by means of suitable conditional random quantities with values in the interval $[0,1]$. 
We  examine the iterated conditional $(B|K)|(A|H)$, by showing that $A|H$ p-entails $B|K$ if and only if $(B|K)|(A|H) = 1$. Then, we show that a p-consistent family $\F=\{E_1|H_1,E_2|H_2\}$ p-entails
 a conditional event $E_3|H_3$ if and only if  $E_3|H_3=1$, or  $(E_3|H_3)|QC(\S)=1$ for some nonempty subset $\S$ of $\F$, where $QC(\S)$ is the quasi conjunction of the conditional events in $\S$. Then, we examine the inference rules \emph{And}, \emph{Cut}, \emph{Cautious Monotonicity}, and  \emph{Or} of System~P  and other well known inference rules (\emph{Modus Ponens}, \emph{Modus Tollens}, \emph{Bayes}).  We also show that $QC(\F)|\C(\F)=1$,  where  $\C(\F)$ is the conjunction of the conditional events in $\F$.
We characterize p-entailment by showing that $\F$ p-entails
 $E_3|H_3$   if and only if $(E_3|H_3)|\C(\F)=1$. Finally, we examine \emph{Denial of the antecedent} and  \emph{Affirmation of the consequent},  where the p-entailment of $(E_3|H_3)$ from $\F$ does not hold, by showing that  $(E_3|H_3)|\C(\F)\neq1.$  
\end{abstract}
\section{Introduction}
The new paradigm psychology of reasoning is characterized by using probability
theory instead of classical bivalent logic as a normative background
theory \cite<see,
e.g.,>{gilio12,oaksford07,over09,elquayam12,pfeifer12y,pfeifer13b,politzer15}. One
of the key topics of the new paradigm psychology of reasoning is how people
interpret and reason about \emph{conditionals} \cite<see, e.g.,>{douven16,edgington95,politzer10,evans04,pfeifer,pfeifer10a,pfeifertulkki17,oaksford03,OverCruz17}. How people interpret
and reason about conditionals was  also one of the key
topics in the (old) logic-based paradigm psychology of reasoning,
which dominated the 
20\textsuperscript{th} century experimental psychology of reasoning. While human interpretation of conditionals
was labeled as ``irrational'' or ``defective'', since the participants'
responses deviated from the semantics of the material conditional,
rationality was revisited and rehabilitated within the new
probabilistic paradigm: specifically, the majority of participants
\begin{itemize}
\item treat
negated antecedents as irrelevant for evaluating whether a 
conditional holds, and 
\item
evaluate their degrees of belief in conditionals by respective conditional
probabilities (and not by the probability of the material conditional). 
\end{itemize}
These findings speak for the conditional event interpretation, and against the material conditional
interpretation, of conditionals. 

 Among
various interpretations of probability, we advocate and use the \emph{coherence-based approach} to probability \cite<see, e.g.,>{BeMR17,biazzo00,biazzo05,CaLS07,coletti02,coletti16,gilio16,GiSa13IJAR,gilio13ins,GiSa14,2012:2UAI,WaPV04},  
which traces back to Bruno \citeA{definetti37,definetti74}. From a
psychological point of view, it is evident that probability serves to measure \emph{degrees of belief} and not some objective quantity in
the  world: this is in line with de Finetti provocative ontological motto
``\emph{Probability
does not exist}''\citeyear[Preface]{definetti74}.  The probabilistic approach based on 
coherence  is thus characterized by \emph{subjective},
 and not by objective, probabilities. Methodologically,
the  approach based on coherence principle differs in many respects to standard
approaches to probabilities. We mention two of them which 
highlight the psychological plausibility of our approach.

 First, contrary to many approaches to probability, the coherence-based approach does not require a complete algebra. 
For drawing a
 probabilistic modus ponens inference, for example, an algebra could be constructed
from the constituents derived from the involved events
in the inference rule.  This is
psychologically plausible, as the reasoning person may focus on only
what is considered to be relevant for drawing the inference. 

Second, conditional probability is a \emph{primitive} notion and it is  not
defined by the fraction of the joint and the marginal probabilities:
the standard definition of $P(C|A)$ by $\frac{P(A \wedge C)}{P(A)}$ requires
to assume 
that $P(A)>0$, as a fraction over zero is undefined. Probabilistic approaches which define conditional probabilities in this
way can therefore not properly manage zero antecedent
probabilities. The subjective probabilistic approach allows for managing zero
antecedent probabilities; moreover, zero probabilities are even exploited for reducing the complexity of the
probabilistic inference. Another aspect of defining conditional
probability \emph{directly} is that the degree of belief in a
conditional \emph{If $A$, then $C$} can be given in a direct way by the reasoner without
presupposing knowledge about $P(A\wedge C)$ and $P(A)$: even as in everyday life it may be impracticable to evaluate the latter two
probabilities,  people
do assess conditionals. For example, if
we want to assess our degree of belief in the conditional that
\emph{If I take the train at six, I am at home at seven}, we can do that
directly, without thinking first about the unconditional probabilities of
\emph{I take the train at six and I am at home at seven} and of
\emph{I take the train at six}. 

In some recent papers of Gilio and Sanfilippo the notions of conjoined and iterated conditionals have been introduced as suitable conditional random quantities \cite{GiSa13c,GiSa13a,GiSa14,GiSa17}. These new objects extend the usual notions of conjunction and conditioning from the case of unconditional events to the case of conditional events. For instance, we developed a semantics for examples like the following \cite<which was presented by>[p. 45]{douven16}: 
\begin{quote}
(I) {\em If the mother is angry if the son gets a B, then she will be furious if the son gets a C}, 
\end{quote}
which is an \emph{iterated (or nested) conditional}. It consists of a conditional in its antecedent\begin{quote}
(A) {\em if the son gets a B, then the mother is angry},
\end{quote}
and a conditonal in its consequent
\begin{quote}
(C) {\em if the son gets a C, then the mother is furious}.
\end{quote}
Of course, the degree of belief in (I) cannot be something like a conditional probability, as the famous triviality results  by \citeA{lewis76} have shown. Rather,  we conceive iterated conditonals like (I) as conditional random quantities (and not as conditonal events) and measure the degree of belief in such objects by previsions $\prev$ \cite<not by probabilities $P$; >{GiSa14,GOPS16,SPOG18}. We will explain the formal details below. 
Interestingly, when we considered the uncertainty propagation rule for the generalized probabilistic modus ponens \cite{ECSQARU17}, where the degree of beliefs are propagated, for instance, from
``$\text{\em The cup broke if dropped}$'' $(A|H)$, 
and 
``if $\text{\em the cup broke if dropped}$, then $\text{\em the cup was fragile}$ $(C|(A|H))$`` to  
``$\text{\em the cup was fragile}$ $(C)$'', we observed, that the uncertainty propagation rules coincide with those of the non-iterated probabilistic modus ponens (i.e., from $P(A)=x$ and $P(C|A)=y$ infer $xy \leq P(C) \leq xy+1-x$). Likewise, we have shown that the uncertainty propagation rules of the iterated version of Centering coincide with the respective (non-iterated) probability propagation rules \cite{SPOG18}. Thus,  a remarkable aspect of the definitions of nested conditionals in terms of conditonal random quantities preserve some well known classical results. 

The main result  of this paper may be  also related to  an analogue result derived from the {\em deduction theorem}. This theorem implies that if an argument is logically valid (or if the premises logically entail the conclusion), then the argument can be transformed into a \emph{logically true conditional}, s.t., the premises are combined by conjunction and form the antecedent and the conclusion forms the consequent of the resulting condtional, which is then a tautology. For example, the logically valid modus ponens (where
 $A\rightarrow C$ denotes the material conditional $\no{A}\vee C$ and
$\models$ denotes logical entailment),
\[
\{A,  A\rightarrow C\} \models C\, ,
\]
can be transformed by the deduction theorem into the following conditional, which is a tautology (and \emph{vice versa}), that is:
\[
(A \wedge (A\rightarrow C))\rightarrow C=\no{(A \wedge (\no{A}\vee C)} \vee C\, =\,\Omega.
\]
Instead of logical entailment, however, we consider in this paper the probabilistic entailment (p-entailment), as introduced by \citeA{adams75,adams98}. Let  $\C(\F)$ denote the conjunction of the conditional events in a p-consistent family $\F$.  We study, in analogy to the deduction theorem,  whether the claim ``a conditional event $E|H$ is p-entailed by a p-consistent family $\F$ of conditional events'' is equivalent to the claim ``the prevision in the iterated conditional $(E|H)|\C(\F)$ is equal to 1''. We examine many cases related to this aspect; in particular, we examine some inference rules of System P and other well known inference rules. 

We remark that this basic relation, between p-entailment and iterated conditioning, appears in its most elementary form when we consider two not impossible  events $A$ and $B$
in the case where  $A \subseteq B$, that is where $A\wedge \no{B}=\emptyset$. In this case $P(A)\leq P(B)$ and then  $A$ p-entails $B$, that is $P(A)=1$ implies $P(B)=1$, and the unique coherent assessment on $B|A$ is $P(B|A)=1$. Therefore, by recalling that in the framework of the betting scheme, when we pay $P(B|A)=x$, we receive $B|A = AB + x\no{A}$; when $A \subseteq B$, it holds that $A$ p-entails $B$ and $B|A = AB + 1 \cdot \no{A} = A + \no{A} = 1$. Conversely, if $B|A=1$, then $P(B|A)=1$; moreover,  
\[
P(B)=P(B|A)P(A)+P(B|\no{A})P(\no{A}) = P(A)+P(B|\no{A})P(\no{A}) \,,
\]
and when $P(A)=1$ it follows that $P(B)=1$, so that $A$ p-entails $B$.

The outline of the paper is as follows.
In Section \ref{SECT:PREL} we give some preliminaries on the notions of coherence and p-entailment for conditional random quantities, which assume values in $[0, 1]$. In Section \ref{SECT:CONG}, after recalling the notions of conjoined  and iterated conditionals, we show that a conditional event $A|H$ p-entails another conditional event $B|K$ if and only if $(B|K)|(A|H) = 1$. Moreover, we show that a p-consistent family of two conditional events $\{E_1|H_1, E_2|H_2\}$ p-entails a conditional event $E_3|H_3$ if and only if  it holds that $(E_3|H_3)|QC(E_1|H_1, E_2|H_2)=1$, where  $QC(E_1|H_1, E_2|H_2)$ denotes the quasi conjunction of $E_1|H_1, E_2|H_2$. We also characterize p-entailment of $E_3|H_3$ from the family $\{E_1|H_1,E_2|H_2\}$ by the property that $E_3|H_3=1$, or
$(E_3|H_3)|QC(\S)=1$ for some nonempty  $\S\subseteq \{E_1|H_1,E_2|H_2\}$.
In Section \ref{SECT:OR} we suitably generalize the notion of iterated conditioning; then, we examine some inference rules of System P and other well known inference rules. The generalization of the notion of iterated conditioning is necessary in order to examine the OR rule. In Section \ref{SECT:MAIN}  we give two results which relate p-entailment and iterated conditioning.  The first result shows that the iterated conditional having as antecedent and consequent 
 the conjunction and the quasi conjunction of two conditional events, respectively,  is equal to 1. The second result characterizes the p-entailment of the conditional event $E_3|H_3$ from a p-consistent family $\{E_1|H_1,E_2|H_2\}$ by the property that the iterated conditional $(E_3|H_3)|((E_1|H_1)\wedge(E_2|H_2)   )$ is equal to 1.
  Finally, we examine two examples  where the  p-entailment of the conditional event $E_3|H_3$ from a p-consistent family $\{E_1|H_1,E_2|H_2\}$ does not hold. We also show that in these cases that  $(E_3|H_3)|((E_1|H_1)\wedge(E_2|H_2))$ does not coincide with 1.
\section{Preliminaries}\label{SECT:PREL}

In our approach  events represent uncertain facts
described by  (non ambiguous) logical propositions.   An event $A$ is
a two-valued logical entity which is either  true ($T$), or false ($F$).
The indicator of an event $A$ is a two-valued numerical quantity which is 1,
or 0, according to  whether $A$ is true, or false, respectively. We use
 the same symbol to refer to an  event and its indicator.
We  denote by
$\Omega$ the sure event and by $\bot$ the impossible one (notice that, when necessary,
the symbol $\bot$  will denote the empty set).
Given two events $A$ and $B$,  we denote by $A\land B$, or simply by $AB$, the  intersection, or conjunction, of $A$ and $B$, as defined in propositional logic; likewise, we denote by  $A \vee B$ the   union, or disjunction, of $A$ and $B$. 
We denote by $\no{A}$ the negation of $A$. 
Of course, the truth values for conjunctions, disjunctions and
negations are defined as usual. 
Given any events $A$ and $B$, we simply write $A \subseteq B$ to denote
that $A$ logically implies $B$, that is  $A\no{B}=\bot$, which
means that it is necessary that $A$ and $\no{B}$ cannot both be true.
	Given two events $A,H$, with $H \neq \emptyset$, the conditional event
	$A|H$ is defined as a three-valued logical entity which is true (T), or
	false (F), or void (V), according to whether $AH$ is true, or
	$\no{A}H$ is true, or $\no{H}$ is true, respectively. 
Given a conditional event $A|H$
with  $P(A|H) = x$, then for (the indicator of) $A|H$ we have
		$A|H = AH + x\no{H} \in \{1,0,x\}$
		\cite[Appendix A.3]{SPOG18}.
We recall below the notion of  logical implication of \citeA{GoNg88} for conditional events   \cite<see also>{gilio13ins}.
\begin{definition}
Given two conditional events $A|H$ and $B|K$ we define that \emph{$A|H$ logically implies $B|K$} (denoted by $A|H \subseteq B|K$) if and only if
 		$AH$ is \emph{true} implies $BK$ is \emph{true} and
 		$\no{B}K$ is \emph{true} implies $\no{A}H$ is 
 		\emph{true}; i.e., $AH\subseteq BK$ and $\no{B}K\subseteq \no{A}H$.
\end{definition}
A generalization of the Goodman and Nguyen logical implication to conditional random quantities has been given by \cite{PeVi14}.\\
The notions of p-consistency and p-entailment of Adams
	\citeyear{adams75} were formulated for conditional events in
	the setting of coherence by  \citeA{gilio10}  \cite<see also>{gilio12a,gilio11ecsqaru,GiSa13IJAR}. 
	\begin{definition}
		\label{PC}
		Let $\mathcal{F}_{n} = \{E_{i}|H_{i} \, , \; i=1,\ldots ,n\}$ be a
		family of $n$ conditional events. Then, $\mathcal{F}_{n}$ is \emph{p-consistent}
		if and only if the probability assessment $(p_{1},p_{2},\ldots ,
		p_{n})=(1,1,\ldots ,1)$ on $\mathcal{F}_{n}$ is coherent.
	\end{definition}
	\begin{definition}
		\label{PE}
		A p-consistent family $\mathcal{F}_{n} = \{E_{i}|H_{i} \, , \; i=1,\ldots ,n\}$
		 \emph{p-entails} a conditional event $E|H$  (denoted by $\mathcal{F}
		_{n} \; \models_{p} \; E|H$) if and only if for any coherent probability
		assessment $(p_{1},\ldots ,
				p_{n},z)$ on $\mathcal{F}_{n} \cup \{E|H
		\}$ it holds that: if $p_{1}=\cdots =p_{n}=1$, then $z=1$.
	\end{definition}
	
	Of course, when $\mathcal{F}_{n}$ p-entails $E|H$, there may be coherent
	assessments $(p_{1},\ldots ,p_{n},z)$ with $z \neq 1$, but in such
	cases $p_{i} \neq 1$ for at least one index $i$. We say that the
	inference from a p-consistent family $\mathcal{F}_{n}$ to $E|H$ is \emph{p-valid} if and only
	if 
	$\mathcal{F}_{n}$ p-entails $E|H$.
We recall  the well known notion of quasi conjunction among conditional events:
\begin{definition}
Given a family  $\mathcal{F}_{n} = \{E_{i}|H_{i} \, , \; i=1,\ldots ,n\}$ of $n$ conditional events, 
the \emph{quasi conjunction} of the conditional events in $\mathcal{F}_{n}$ is defined as
\[
QC(\mathcal{F}_n)= \bigwedge_{i=1}^n(\no{H}_i\vee E_iH_i)| (\bigvee_{i=1}^n H_i).
\]	
\end{definition}

Moreover, we recall the following characterization of  p-entailment \cite{GiSa13IJAR}:
\begin{theorem}\label{ENTAIL-CS}{\rm
Let  a p-consistent family $\mathcal{F}_n = \{E_1|H_1, \ldots, E_n|H_n\}$ and a conditional event $E|H$ be given. The following assertions are equivalent: \\
1. $\mathcal{F}_n$ p-entails $E|H$; \\
 2. The assessment $\mathcal{P}=(1,\ldots,1,z)$ on $\mathcal{F}=\mathcal{F}_n \cup \{E|H\}$, where $P(E_i|H_i)=1,$ $i=1,\ldots,n, P(E|H)=z$, is coherent if and only if $z=1$; \\
 3. The assessment $\mathcal{P}=(1,\ldots,1,0)$ on $\mathcal{F}=\mathcal{F}_n \cup \{E|H\}$, where $P(E_i|H_i)=1,$ $i=1,\ldots,n, P(E|H)=0$, is not coherent; \\
 4. Either there exists a nonempty $\mathcal{S} \subseteq \mathcal{F}_n$  such that $QC(\mathcal{S})$ implies $E|H$, or $H \subseteq E$; \\
5.  There exists a nonempty $\mathcal{S} \subseteq \mathcal{F}_n$ such that $QC(\mathcal{S})$ p-entails $E|H$.
}\end{theorem}
We also recall the characterization of the p-entailment for two conditional events \cite[Theorem 7]{gilio13ins}:
\begin{theorem}\label{THM:pentailmentfor2}
Given two conditional events $A|H$, $B|K$, with $AH\neq \emptyset$. It holds that
\[
A|H \Rightarrow_p B|K \Longleftrightarrow A|H\subseteq B|K, \mbox{  or } K\subseteq B
\Longleftrightarrow \Pi\subseteq \{(x,y)\in[0,1]^2: x\leq y\},
\]
where $\Pi$ is the set of coherent assessments $(x,y)$ on $\{A|H,B|K\}$.
\end{theorem}
We denote by $X$ a random quantity, that is  an 
uncertain real quantity,  which has a well determined but unknown value. 
We assume that  $X$ has a finite set of possible values. Given any event $H\neq \emptyset$, 
agreeing to the betting metaphor, if you  assess that the prevision of $``X$ {\em conditional on} $H$'' (or short:  $``X$ {\em given} $H$''), $\pr(X|H)$, is equal to $\mu$, this means that for any given  real number $s$ you are willing to pay an amount $\mu s$ and to receive  $sX$, or $\mu s$, according  to whether $H$ is true, or  false (i.e., when the bet is called off), respectively. 
In particular, when $X$ is (the indicator of) an event $A$, then $\prev(X|H)=P(A|H)$. 
In \citeA{GiSa14} the notions of conjoined, disjoined,  and iterated conditionals  have been studied in the framework of conditional  random quantities. In particular, the next result establishes  some conditions under which two conditional random quantities $X|H$ and  $Y|K$ coincide \cite[Theorem 4]{GiSa14}:
\begin{theorem}\label{EQ-CRQ}{\rm Given any events $H\neq \emptyset$ and  $K\neq \emptyset$, and any random quantities $X$ and $Y$, let $\Pi$ be the set of the coherent prevision assessments $\pr(X|H)=\mu$ and $\pr(Y|K)=\nu$. \\
(i) Assume that, for every $(\mu,\nu)\in \Pi$, the values of  $X|H$ and $Y|K$ always  coincide when  $H\vee K$ is true; then   $\mu=\nu$ for every $(\mu,\nu)\in \Pi$. \\
(ii) For every $(\mu,\nu)\in \Pi$,  the values of  $X|H$ and $Y|K$ always coincide when  $H\vee K$ is true  if and only if $X|H=Y|K$.
}\end{theorem}

\section{Generalized System P and Compound Conditionals}
\label{SECT:CONG}
In this  section we recall the notions of conjunction  and iterated conditioning.  Then, we show that  $A|H$ p-entails  $B|K$ if and only if $(B|K)|(A|H) = 1$. Moreover, we show  that $\{E_1|H_1, E_2|H_2\}$ p-entails $E_3|H_3$ if and only if $(E_3|H_3)|QC(E_1|H_1, E_2|H_2)=1$.
\subsection{Exploring  conjunction  and iterated conditioning}
We recall below the definition of conjuntion of two conditional events $A|H$ and $B|K$ \cite{GiSa13a,GiSa13c,GiSa14}.
Different approaches to compounded conditionals, not based on coherence,  have been developed by other authors  \cite<see, e.g.,>{Kauf09,mcgee89}.

\begin{definition}\label{CONJUNCTION}{\rm Given any pair of conditional events $A|H$ and $B|K$, with $P(A|H)=x$ and $P(B|K)=y$, we define their conjunction
		as the conditional random quantity $(A|H) \wedge (B|K) = Z\,|\, (H \vee K)$, where $Z=\min \, \{A|H, B|K\}$.
}\end{definition}
In betting terms, $z=\mathbb{P}[(A|H)\wedge(B|K)]$ represents the amount you agree to pay, with the proviso that you will receive the quantity:
\begin{equation}\label{EQ:CONJUNCTION}
(A|H)\wedge(B|K) =\left\{\begin{array}{ll}
1, &\mbox{if $AHBK$ is true,}\\
0, &\mbox{if $\no{A}H\vee \no{B}K$ is true,}\\
x, &\mbox{if $\no{H}BK$ is true,}\\
y, &\mbox{if $AH\no{K}$ is true,}\\
z, &\mbox{if $\no{H}\no{K}$ is true}.
\end{array}
\right.
\end{equation}
From (\ref{EQ:CONJUNCTION}), it follows that the conjunction  $(A|H) \wedge (B|K)$ is the following  random quantity
\begin{equation}\label{EQ:REPRES}
(A|H) \wedge (B|K)=1 \cdot AHBK + x \cdot \no{H}BK + y \cdot AH\no{K} + z \cdot \no{H}\no{K}\,.
\end{equation} 
We observe that if $H=K$, 
then $\no{H}BK= AH\no{K}=\emptyset$, so that
$(A|H) \wedge (B|K)=ABH + z\no{H};$
moreover, $AB|H=ABH+p\no{H}$, where $p=P(AB|H)$. 
We notice that $(A|H)\wedge (B|H)$ and $AB|H$ coincide when $H$ is true; then,
by Theorem \ref{EQ-CRQ}, $z=p$; thus,  
\begin{equation}\label{EQ:H=K}
(A|H)\wedge (B|H)=AB|H.
\end{equation}
We recall that, given any coherent assessment $(x,y)$ on $\{A|H, B|K\}$, with $A,H,B,K$ logically independent, and with $H \neq \bot, K \neq \bot$, the extension $z = \mathbb{P}[(A|H) \wedge (B|K)]$ is coherent if and only if  the following Fr\'echet-Hoeffding bounds are satisfied \cite[Theorem 7]{GiSa14}:
		\begin{equation}\label{LOW-UPPER}
		max\{x+y-1,0\} = z' \; \leq \; z \; \leq \; z'' = min\{x,y\} \,.
		\end{equation}
Note that the bounds in~(\ref{LOW-UPPER}) coincide with the bounds for the conjunction of  unconditional probabilities (i.e.,  {\em if $P(A)=x$ and $P(B)=y$, then $max\{x+y-1,0\} \leq P(AB) \leq min\{x,y\}$}).

	We now turn to recalling and discussing the notion of iterated conditioning (see, e.g., \cite{GiSa13c,GiSa13a,GiSa14}).
\begin{definition}[Iterated conditioning]
	\label{DEF:ITER-COND} Given any pair of conditional events $A|H$ and $B|K$, with $AH\neq \emptyset$, the iterated
	conditional $(B|K)|(A|H)$ is defined as the conditional random quantity
\begin{equation}\label{EQ:ITER-COND}
	(B|K)|(A|H) = (B|K) \wedge (A|H) + \mu \no{A}|H,
\end{equation}
	where
	$\mu =\mathbb{P}[(B|K)|(A|H)]$.
\end{definition}
\begin{remark}\label{REM:A|H=0}
Notice that we assumed that $AH\neq \emptyset$ to give a nontrivial meaning to the notion of the iterated conditional. Indeed,   if  $AH$ were equal to $\emptyset$,  that is $A|H=0$, then it would be the case that $\no{A}|H =1$ and     
$(B|K)|(A|H)=(B|K)|0=(B|K) \wedge (A|H) + \mu \no{A}|H=\mu$ would follow; that is, $(B|K)|(A|H)$ 
would coincide with the (indeterminate) value $\mu$. Similarly in the case of  $B|\emptyset$   (which is of no interest): the trivial iterated conditional $(B|K)|0$
is not considered in our approach.
\end{remark}
We observe that, by linearity of prevision, it holds that 
\[
\mu = \prev((B|K)|(A|H)) = \prev((B|K) \wedge (A|H)) + \mu P(\no{A}|H) = z + \mu(1-x) \,,
\]
from which it follows that $z=\mu x$. Here, when $x>0$, we obtain $\mu = \frac{z}{x} \in [0,1]$. Notice that $z + \mu(1-x)$, i.e. $\mu$, is the value of $(B|K)|(A|H)$ when $\no{H}\no{K}$ is true. Then, by observing that 
\[
\no{A}H\no{K} \vee \no{A}HBK \vee \no{A}H\no{B}K \vee \no{H}\no{K} = \no{A}H \vee \no{H}\no{K} \,,
\]
we obtain
\[
(B|K)|(A|H) = \left\{\begin{array}{ll}
1, &\mbox{if $AHBK$ is true,}\\
0, &\mbox{if $AH\no{B}K$ is true,}\\
y, &\mbox{if $AH\no{K}$ is true,}\\
x + \mu(1-x), &\mbox{if $\no{H}BK$ is true,}\\
\mu(1-x), &\mbox{if $\no{H}\no{B}K$ is true,}\\
\mu,  &\mbox{if $\no{A}H\no{K}$ is true,}\\
\mu,  &\mbox{if $\no{A}HBK$ is true,}\\
\mu,  &\mbox{if $\no{A}H\no{B}K$ is true,}\\
\mu, &\mbox{if $\no{H}\no{K}$ is true}.
\end{array}
\right.
= \left\{\begin{array}{ll}
1, &\mbox{if $AHBK$ is true,}\\
0, &\mbox{if $AH\no{B}K$ is true,}\\
y, &\mbox{if $AH\no{K}$ is true,}\\
x + \mu(1-x), &\mbox{if $\no{H}BK$ is true,}\\
\mu(1-x), &\mbox{if $\no{H}\no{B}K$ is true,}\\
\mu,  &\mbox{if $\no{A}H \vee \no{H}\no{K}$ is true.}
\end{array}
\right.
%
%
%
%
\]

In particular, when $x=0$, it holds that
\[
(B|K)|(A|H) = \left\{\begin{array}{ll}
1, &\mbox{if $AHBK$ is true,}\\
0, &\mbox{if $AH\no{B}K$ is true,}\\
y, &\mbox{if $AH\no{K}$ is true,}\\
\mu, &\mbox{if $\no{H}BK$ is true,}\\
\mu, &\mbox{if $\no{H}\no{B}K$ is true,}\\
\mu,  &\mbox{if $\no{A}H \vee \no{H}\no{K}$ is true.}
\end{array}
\right.
= \left\{\begin{array}{ll}
1, &\mbox{if $AHBK$ is true,}\\
0, &\mbox{if $AH\no{B}K$ is true,}\\
y, &\mbox{if $AH\no{K}$ is true,}\\
\mu, &\mbox{if $\no{A}H \vee \no{H}$ is true.}
\end{array}
\right.
\]
As we can see, in order that the prevision assessment  $\mu$ on $(B|K)|(A|H)$ be coherent, $\mu$ must belong to the convex hull of the values $0,y,1$; that is, (also when $x=0$) it must be that $\mu \in [0,1]$. 
\begin{proposition}\label{PROP:ANDGN}
Given two conditional events $A|H$ and $B|K$, it holds that
\begin{equation}
A|H \subseteq B|K \Longrightarrow (A|H)\wedge (B|K)=A|H\,.
\end{equation}
\end{proposition}
\begin{proof}
We set $P(A|H)=x$, $P(B|K)=y$, and $\prev[(A|H)\wedge (B|K)]=z$.
As $A|H \subseteq B|K$, it holds that  $AH\no{B}K=AH\no{K}=\no{H}\no{B}K=\emptyset$ and $AHBK=AH$
\cite[Remark 3]{gilio13ins}. Then,
\[
(A|H)\wedge (B|K)=AHBK+x \no{H}BK+y\no{K}AH+z\no{H}\no{K}=AH+x \no{H}BK+z\no{H}\no{K}.
\]
Moreover,
\[
A|H=AH+x\no{H}=AH+x\no{H}BK+x\no{H}\no{K}.
\]
We notice that $(A|H)\wedge (B|K)$ and $A|H$ coincide when $H\vee K$ is true. Then,
  $z=x$ follows  from Theorem~\ref{EQ-CRQ}. Therefore,  $(A|H)\wedge (B|K)=A|H$.
\qed
\end{proof}

The following theorem shows that a  conditional $A|H$ p-entails another conditional $B|K$  
if and only if the unique coherent prevision assessment for the corresponding iterated conditonal $(B|K)|(A|H)$ is equal to one.
\begin{theorem}\label{THM:p-entail-iter}
Given two (p-consistent) conditional events $A|H$ and $B|K$, it holds that,
\begin{equation}
A|H \Rightarrow_p B|K \Longleftrightarrow (B|K)|(A|H)=1.
\end{equation}
\end{theorem}
\begin{proof}
$(\Rightarrow)$.
We distinguish two cases: $(i)$ $A|H \subseteq B|K$; $(ii)$ $K\subseteq B$.
Case $(i)$.  We remark that if $A|H \subseteq B|K$, then $A|H \leq B|K$ and $P(A|H) \leq P(B|K)$; moreover, $(A|H) \wedge (B|K) = A|H$. Then, by defining $\prev((B|K)|(A|H)) = \mu,\, P(A|H)=x$, we obtain
\[
(B|K)|(A|H) = (A|H) \wedge (B|K) + \mu \no{A}|H = A|H + \mu \no{A}|H  = \left\{\begin{array}{ll}
1, &\mbox{if $AH$ is true,}\\
\mu, &\mbox{if $\no{A}H$ is true,} \\
x+\mu(1-x), &\mbox{if $\no{H}$ is true.}
\end{array}
\right.
\]
By linearity of prevision, we obtain
\begin{equation} \label{EQ:PREVITER}
\prev((B|K)|(A|H)) = \mu = P(A|H) + \mu P(\no{A}|H) = x + \mu(1-x) \,;
\end{equation}
which implies that
\[
(B|K)|(A|H) = \left\{\begin{array}{ll}
1, &\mbox{if $AH$ is true,}\\
\mu, &\mbox{if $\no{A}H \vee \no{H}$ is true.}
\end{array}
\right.
\]
In order for $\mu$ to be coherent, $\mu$ must belong to the convex hull of the set $\{1\}$; i.e. $\mu=1$. In other words, given two conditional events $A|H$ and $B|K$, with $A|H \subseteq B|K$, it holds that: $\prev((B|K)|(A|H))=1$.  Thus $(B|K)|(A|H)=1$. \\
Case $(ii)$. If $K\subseteq B$ it holds that $P(B|K)=y=1$ and $B|K=1$.
Then, 
  $(A|H)\wedge (B|K)=(A|H)|(H\vee K)=A|H$ \cite<see>[Remark 4]{GiSa13c}.
Moreover,   $(B|K)|(A|H)=A|H+\mu\no{A}|H$ and by linearity of prevision it holds that
$\mu=x+\mu(1-x)$. Then, 
\[
(B|K)|(A|H) = \left\{\begin{array}{ll}
1, &\mbox{if $AH$ is true,}\\
\mu, &\mbox{if $\no{A}H$ is true,} \\
x+\mu(1-x), &\mbox{if $\no{H}$ is true.}
\end{array}
\right.
=
\left\{\begin{array}{ll}
1, &\mbox{if $AH$ is true,}\\
\mu, &\mbox{if $\no{A}H \vee \no{H}$ is true.}
\end{array}
\right.
\]
Then, by coherence, $\mu=1$ and  $(B|K)|(A|H)=1$.\\
Thus,  p-entailment of $B|K$ from   $A|H$ implies
 $(B|K)|(A|H)=1$.
  
 $(\Leftarrow)$. Assume that $(B|K)|(A|H)=1$, so that the unique coherent assessment for  $\prev[(B|K)|(A|H)]$ is $\mu=1$.
 Then, by observing that $ \prev[(A|H)\wedge (B|K)]\leq P(B|K)=y$, it follows that
 \[
 \prev[(A|H)\wedge (B|K)]= \prev[ (B|K)|(A|H)]P(A|H)=P(A|H)=x\leq y.
 \] 
 Then, when $x=1$, it holds that $y=1$; that is, $A|H$ p-entails $B|K$. 
\qed
\end{proof}	
\begin{corollary}\label{COR:QC}
	Let  three conditional events $E_1|H_1$, $E_2|H_2$, and  $E_3|H_3$  be given, where $\{E_1|H_1, E_2|H_2\}$ is p-consistent.
	The quasi conjunction $QC(E_1|H_1, E_2|H_2)$ p-entails $E_3|H_3$ if and only if 
	$(E_3|H_3)|QC(E_1|H_1,E_2|H_2) = 1$.
\end{corollary}
\begin{proof} The assertion directly follows by applying Theorem \ref{THM:p-entail-iter}, with $A|H = QC(E_1|H_1,E_2|H_2)$ and $B|K = E_3|H_3$. \qed
\end{proof} 
In the next result we characterize the p-entailment of $E_3|H_3$ from the family $\{E_1|H_1,E_2|H_2\}$ by the property that $E_3|H_3=1$, or
$(E_3|H_3)|QC(\S)=1$ for some nonempty  $\S\subseteq \{E_1|H_1,E_2|H_2\}$.
\begin{theorem}
	Let   three conditional events $E_1|H_1$, $E_2|H_2$, and $E_3|H_3$ be given, where $\{E_1|H_1, E_2|H_2\}$ is p-consistent. Then, the family  $\{E_1|H_1, E_2|H_2\}$ p-entails $E_3|H_3$ if and only if  at least one of the following conditions is satisfied: $(i)$ $E_3|H_3=1$; $(ii)$ $(E_3|H_3)|(E_1|H_1)=1$;
	$(iii)$ $(E_3|H_3)|(E_2|H_2)=1$; 	$(iv)$ 	$(E_3|H_3)|QC(E_1|H_1,E_2|H_2) = 1$.
\end{theorem}
\begin{proof}
$(\Rightarrow)$. By Theorem $\ref{ENTAIL-CS}$, 
as $\{E_1|H_1, E_2|H_2\}$ p-entails $E_3|H_3$, it follows that  $QC(\mathcal{S})\subseteq E_3|H_3$ for some $\emptyset \neq  \mathcal{S}\subseteq \{E_1|H_1, E_2|H_2\}$, or $H_3 \subseteq E_3$.
If $H_3 \subseteq E_3$, then $P(E_3|H_3)=1$ and $E_3|H_3=1$.
If $\S=\{E_i|H_i\}$, for $i=1$ or $i=2$, by Theorem~\ref{THM:p-entail-iter} it holds that $(E_3|H_3)|(E_i|H_i)=1$.
If $\S=\{E_1|H_1, E_2|H_2\}$, then by Corollary~\ref{COR:QC} it holds that 	$(E_3|H_3)|QC(E_1|H_1,E_2|H_2) = 1$.\\
$(\Leftarrow)$. 
If $E_3|H_3=1$ then the unique coherent assessment on $E_3|H_3$ is $P(E_3|H_3)=1$. This means that $H_3\subseteq E_3$ and then $\{E_1|H_1, E_2|H_2\}$ p-entails $E_3|H_3$.\\
If $(E_3|H_3)|(E_i|H_i)=1$, for $i=1$ or $i=2$, then by Theorem~\ref{THM:p-entail-iter} it holds that 
$E_i|H_i$ p-entails $E_3|H_3$ and hence,  by Theorem~\ref{THM:p-entail-iter},  $\{E_1|H_1, E_2|H_2\}$ p-entails $E_3|H_3$.\\
Finally, if $(E_3|H_3)|QC(E_1|H_1,E_2|H_2) = 1$, then by Corollary~\ref{COR:QC} it holds that 
$QC(E_1|H_1,E_2|H_2)$ p-entails $E_3|H_3$ and hence,  by Theorem~\ref{THM:p-entail-iter},  $\{E_1|H_1, E_2|H_2\}$ p-entails $E_3|H_3$.
\qed
\end{proof}

\section{Iterated conditionals and some inference rules}\label{SECT:RULES}
\label{SECT:OR}
In this section we  examine some  inference rules with  $\{E_1|H_1, E_2|H_2\}$ as the premise set,
 and $E_3|H_3$  as the  conclusion, by showing that, if  $\{E_1|H_1, E_2|H_2\}\Rightarrow_p  E_3|H_3$, then  $(E_3|H_3)|((E_1|H_1)\wedge (E_2|H_2))=1$. 
  The notion of conjunction of three conditional events is given below \cite{GiSa17}.
 \begin{definition}
 	Given a family of three conditional events $\F=\{E_1|H_1$, $E_2|H_2$,$E_3|H_3$\}, we set $P(E_i|H_i) = x_i$, $i=1,2,3$, $\mathbb{P}[(E_i|H_i)\wedge(E_j|H_j)]=x_{ij}=x_{ji}$, $i\neq j$.
 	The conjunction  $\C(\F)=(E_1|H_1) \wedge (E_2|H_2)\wedge (E_3|H_3)$ is defined as the conditional random quantity
 	\begin{equation}\label{EQ:CONJUNCTION3}
\C(\F)=(E_1|H_1) \wedge (E_2|H_2)\wedge (E_3|H_3)=	
 	\left\{
 	\begin{array}{llll}
 	1, &\mbox{ if } E_1H_1E_2H_2E_3H_3 \mbox{ is true}\\
 	0, &\mbox{ if } \no{E}_1H_1 \vee \no{E}_2H_2 \vee \no{E}_3H_3 \mbox{ is true},\\
 	x_1,& \mbox{ if } \no{H}_1E_2H_2E_3H_3 \mbox{ is true},\\
 	x_2,& \mbox{ if } \no{H}_2E_1H_1E_3H_3 \mbox{ is true},\\
 	x_3, &\mbox{ if } \no{H}_3E_1H_1E_2H_2 \mbox{ is true}, \\
 	x_{12}, &\mbox{ if } \no{H}_1\no{H}_2E_3H_3 \mbox{ is true}, \\
 	x_{13}, &\mbox{ if } \no{H}_1\no{H}_3E_2H_2 \mbox{ is true}, \\
 	x_{23}, &\mbox{ if } \no{H}_2\no{H}_3E_1H_1 \mbox{ is true}, \\
 	x_{123}, &\mbox{ if } \no{H}_1\no{H}_2\no{H}_3 \mbox{ is true} \\
 	\end{array}
 	\right.
 	\end{equation}
 	where   $x_{123}=\mathbb{P}[\C(\F)]$.	
 	\end{definition}
 We recall below the definition of the object $(E_3|H_3)|((E_1|H_1)\wedge (E_2|H_2))$, which is under study in \cite{GiSa17wp}.
\begin{definition}\label{DEF:GENITER}
	Let be given three conditional events $E_1|H_1$, $E_2|H_2,$ and $E_{3}|H_{3}$, with $(E_1|H_1)  \wedge (E_2|H_2)\neq 0$. We denote by $(E_{3}|H_{3})|((E_1|H_1)  \wedge (E_2|H_2))$ the conditional random quantity
	\[
	(E_1|H_1) \wedge (E_2|H_2)\wedge (E_{3}|H_{3}) + \mu (1-(E_1|H_1) \wedge  (E_2|H_2)) \,,
	\]
	where $\mu = \prev[(E_{3}|H_{3})|((E_1|H_1) \wedge  (E_2|H_2))]$.
\end{definition}
\begin{remark}\label{REM:COMPOUND}
	We observe that, defining 	$\prev[(E_1|H_1) \wedge (E_2|H_2)\wedge (E_{3}|H_{3})]=t$ and 
	$\prev[(E_1|H_1) \wedge (E_2|H_2)]=z$, by the linearity of prevision it holds that $\mu=t+\mu(1-z)$; then, $t=\mu z$, that is
	\[
	\prev[(E_1|H_1) \wedge (E_2|H_2)\wedge (E_{3}|H_{3})]=
	\prev[(E_{3}|H_{3})|((E_1|H_1) \wedge (E_2|H_2))]\prev[(E_1|H_1) \wedge (E_2|H_2)].
	\]
\end{remark}

\paragraph{Modus Ponens:}
$\{C|A, A\}\Rightarrow_p C$.
It holds that $(C|A)\wedge A= AC= QC(AC)$; then,
by Theorem \ref{THM:p-entail-iter},  as $AC\subseteq C$   it follows that
\[
C| ((C|A)\wedge A)=
C| (QC((C|A),A)=
C|AC=  1\,.
\]
This can be seen as an analogy to the  fact that the modus ponens is logically valid in logic and that the   probabilistic modus ponens is p-valid.
\paragraph{Modus Tollens:}
$\{C|A, \no{C}\}\Rightarrow_p \no{A}$.
It holds that $(C|A)\wedge \no{C}= x\no{A}\no{C}$,  where $x=P(C|A)$, while  $QC(C|A,\no{C})=\no{A}\no{C}$; then, assuming $x>0$, we obtain
\[
\no{A}|((C|A)\wedge \no{C}) =\no{A}\wedge (C|A)\wedge \no{C}+\mu(1-(C|A)\wedge \no{C})) = \left\{\begin{array}{ll}
\mu, &\mbox{if $A\vee C$ is true,}\\
x+\mu(1-x), &\mbox{if $\no{A}\no{C}$ is true.}
\end{array}
\right.
\]
By coherence it must be the case that $\mu=x+\mu(1-x)$,  i.e.,    $x=\mu x$, which implies  $\mu=x+\mu(1-x)=1$; therefore, 
\[
\no{A}|((C|A)\wedge \no{C})=1\,.
\]
This can be seen as an analogy to the  fact that the modus tollens is logically valid in logic and that the   probabilistic modus tollens is p-valid. 
Notice that, if $x= 0$, then $((C|A)\wedge \no{C})=0$ and the object $\no{A}|((C|A)\wedge \no{C})=\no{A}|0=\mu$, which is indeterminate (see Remark \ref{REM:A|H=0}). 

\paragraph{Bayes.}
We note that $(E|AH)\wedge (H|A)= EH|A=QC(E|AH,H|A)$; 
then, as $EH|A\subseteq H|EA$, by Theorem (\ref{THM:pentailmentfor2}) it holds that $\{E|AH, H|A\} \Rightarrow_p  H|EA$. Moreover, 
by Theorem \ref{THM:p-entail-iter}, it follows that
\[
(H|EA)| ((E|HA)\wedge H|A)=
(H|EA)| (QC((E|HA),H|A)=
(H|EA)|(EH|A)=  1\,.
\]
In particular, if $A=\Omega$, we obtain $(H|E)|(EH)=  1$.

\subsection{And, Cut, and Cautious Monotonicity of  System P}
In this section we consider the following inference rules of System P \cite{kraus90}: And, Cut, and Cautious Monotonicity (short: CM). System P is a basic nonmonotonic reasoning which allows for retracting conclusions in the light of new premises. The probabilistic versions of the rules of System P are p-valid \cite{adams75,biazzo02,gilio02}. Experimental evidence supports the   psychological plausibility of System P \cite<see, e.g.>{dasilva02,pfeifer03,pfeifer,schurz05}. 
\paragraph{And rule:}
$\{B|A, C|A\}\Rightarrow_p BC|A$.
By formula (\ref{EQ:H=K}),  it holds that $(B|A)\wedge (C|A)= BC|A= QC(B|A,C|A)$; then,
  by Theorem \ref{THM:p-entail-iter},  as $BC|A\subseteq BC|A$   it follows that
\[
 (BC|A)| ((C|A)\wedge (B|A))=
 (BC|A)| QC(B|A,C|A)=
  (BC|A)| (BC|A)=  1\,.
 \]
\paragraph{Cut rule:}
 $\{C|AB, B|A\}\Rightarrow_p C|A$. 
We note that  $(C|AB)\wedge (B|A)=BC|A=QC(C|AB,B|A)$; then,
  by Theorem \ref{THM:p-entail-iter}, as $BC|A \subseteq C|A $ it follows that 
\[
 (C|A)| ((C|AB)\wedge (B|A)) =(C|A)|QC(C|AB,B|A) = (C|A)|(BC|A)=1\,.
 \]
\paragraph{CM rule:}
 $\{C|A, B|A\}\Rightarrow_p C|AB$. 
By formula (\ref{EQ:H=K}),  it holds that $(C|A)\wedge (B|A)=BC|A=QC(C|A,B|A)$; then,
  by Theorem \ref{THM:p-entail-iter}, as $BC|A\subseteq C|AB$ it follows that  
\[
 (C|AB)| ((C|A)\wedge (B|A))= (C|AB)|QC(C|A,B|A)= (C|AB)|(BC|A) = 1\,.
 \]

\subsection{Or rule of System P}
 We recall that the Or rule is p-valid, that is  
$\{C|A, C|B\}\Rightarrow_p C|(A \vee B)$. The next result   shows that  the conclusion of the Or rule, $C|(A \vee B)$,  given the conjunction of the premises, $(C|A) \wedge (C|B)$, coincides with 1.
\begin{theorem}	\label{THM:OR}
	Given a $p$-consistent family $\{C|A, C|B\}$ it holds that 
\[
(C|(A \vee B)|((C|A) \wedge (C|B))=1.
\]
\end{theorem}
\begin{proof}
By Definition \ref{DEF:GENITER}, we obtain
\[
(C|(A \vee B))|((C|A) \wedge (C|B)) = (C|(A \vee B)) \wedge (C|A) \wedge (C|B) + \mu [1-(C|A) \wedge (C|B)] \,,
\]
where $\mu=\prev[(C|(A \vee B))|(C|A) \wedge (C|B)]$. We set $P(C|A)=x$, $P(C|B)=y$, and $\prev((C|A) \wedge (C|B))=z$; then, 
\[
(C|A) \wedge (C|B) = \left\{\begin{array}{ll}
1, &\mbox{if $ABC$ is true,}\\
0, &\mbox{if $(A\vee B)\no{C}$ is true,}\\
x, &\mbox{if $\no{A}BC$ is true,}\\
y, &\mbox{if $A\no{B}C$ is true,}\\
z, &\mbox{if $\no{A}\no{B}$ is true.}
\end{array}
\right.
\]
Moreover, by defining $\prev[(C|(A \vee B)) \wedge (C|A) \wedge (C|B)]=t$, we obtain
\[
(C|(A \vee B)) \wedge (C|A) \wedge (C|B) = \left\{\begin{array}{ll}
1, &\mbox{if $ABC$ is true,}\\
0, &\mbox{if $(A\vee B)\no{C}$ is true,}\\
x, &\mbox{if $\no{A}BC$ is true,}\\
y, &\mbox{if $A\no{B}C$ is true,}\\
t, &\mbox{if $\no{A}\no{B}$ is true.}
\end{array}
\right.
\]
As we can see, $(C|(A \vee B)) \wedge (C|A) \wedge (C|B)$ and $(C|A) \wedge (C|B)$ coincide when $A \vee B$ is true; then, by Theorem \ref{EQ-CRQ}  it holds that $t = z$, so that 
\[
(C|(A \vee B)) \wedge (C|A) \wedge (C|B) = (C|A) \wedge (C|B).
\]
 Then,
\begin{equation} \label{EQ:ORITER}
(C|(A \vee B))|((C|A) \wedge (C|B)) = (C|A) \wedge (C|B) + \mu [1-(C|A) \wedge (C|B)] \,,
\end{equation}
and by the linearity of prevision we obtain $\mu=z+\mu(1-z)$, so that $z=\mu z$. Moreover, by (\ref{EQ:ORITER}) we obtain
\[
\begin{array}{lll}
(C|(A \vee B))|((C|A) \wedge (C|B)) &=& \left\{\begin{array}{ll}
1, &\mbox{if $ABC$ is true,}\\
x+\mu(1-x), &\mbox{if $\no{A}BC$ is true,}\\
y+\mu(1-y), &\mbox{if $A\no{B}C$ is true,}\\
\mu, &\mbox{if $AB\no{C}$ is true,}\\
\mu, &\mbox{if $\no{A}B\no{C}$ is true,}\\
\mu, &\mbox{if $A\no{B}\no{C}$ is true,}\\
\mu, &\mbox{if $\no{A}\no{B}$ is true,}
\end{array}
\right.
\end{array}
\]
which reduces to 
\[
\begin{array}{lll}
(C|(A \vee B))|((C|A) \wedge (C|B))=\left\{\begin{array}{ll}
1, &\mbox{if $ABC$ is true,}\\
x+\mu(1-x), &\mbox{if $\no{A}BC$ is true,} \\
y+\mu(1-y), &\mbox{if $A\no{B}C$ is true,}\\
\mu, &\mbox{if $\no{A}\no{B}C \vee \no{C}$ is true}. 
\end{array}
\right.
\end{array}
\]
In order to prove that $(C|(A \vee B))|((C|A) \wedge (C|B))=1$, we distinguish the following cases: (a) $z>0$; (b) $z=x=y=0$; (c)  $z=0, x>0, y>0$; (d) $z=y=0, x>0$; (e) $z=x=0, y>0$. \\
Case (a). By recalling that $z=\mu z$, as $z>0$ it follows that  $\mu=1$; then, $y+\mu(1-y)=x+\mu(1-x)=1$, so that
\[
(C|(A \vee B))|((C|A) \wedge (C|B)) = \left\{\begin{array}{ll}
1, &\mbox{if $ABC \vee A\no{B}C \vee \no{A}BC$ is true,}\\
\mu, &\mbox{if $\no{A}\no{B}C \vee \no{C}$ is true.}
\end{array}
\right.
\]
Then, by coherence, $\mu = 1$ and $(C|(A \vee B))|((C|A) \wedge (C|B)) = 1$.\\
Case (b). As $x=y=0$, it holds that $x+\mu(1-x)=y+\mu(1-y)=\mu$; then 
\[
(C|(A \vee B))|((C|A) \wedge (C|B)) = \left\{\begin{array}{ll}
1, &\mbox{if $ABC$ is true,}\\
\mu, &\mbox{if $\no{ABC}$ is true.}
\end{array}
\right.\,,
\]
and, by coherence, $\mu=1$; thus, $(C|(A \vee B))|((C|A) \wedge (C|B)) =1$. \\
Case (c). By coherence, $\mu$ is a linear convex combination of the values $1$, $y+\mu(1-y)$, and $x+\mu(1-x)$, that is,
\begin{equation}\label{EQ:CONV}
\mu=\lambda_1+\lambda_2(y+\mu(1-y))+\lambda_3(x+\mu(1-x)) \,,
\end{equation}
with $\lambda_h \geq 0, h=1,2,3,$ and $\lambda_1+\lambda_2+\lambda_3= 1$.
The equation (\ref{EQ:CONV}) can be written as
\[
\mu(\lambda_1+\lambda_2y+\lambda_3x) =\lambda_1+\lambda_2y+\lambda_3x \,,
\]
where $\lambda_1+\lambda_2y+\lambda_3x>0$; then
$\mu=y+\mu(1-y)=x+\mu(1-x)=1$ and $(C|(A \vee B))|((C|A) \wedge (C|B)) = 1$. \\
Case (d). As $y=0$ it holds that $y+\mu(1-y)=\mu$; then,
\[
(C|(A \vee B))|((C|A) \wedge (C|B)) = \left\{\begin{array}{ll}
1, &\mbox{if $ABC$ is true,}\\
x+\mu(1-x), &\mbox{if $\no{A}BC$ is true,} \\
\mu, &\mbox{if $\no{BC}$ is true}.
\end{array}
\right.
\]
By coherence, $\mu$ is a linear convex combination of the values $1, x+\mu(1-x)$, that is
\begin{equation}\label{EQ:CONV1}
\mu=\lambda_1+\lambda_2(x+\mu(1-x)) \,,\;\; \lambda_1 \geq 0 \,,\; \lambda_2 \geq 0 \,,\; \lambda_1 + \lambda_2 = 1 \,.
\end{equation}
The equation (\ref{EQ:CONV1}) can be written as
$\mu(\lambda_1+\lambda_2x) =\lambda_1+\lambda_2x$,
where $\lambda_1+\lambda_2x>0$; then, $\mu=x+\mu(1-x)=1$ and $(C|(A \vee B))|((C|A) \wedge (C|B)) =1$.\\
Case (e). Since $x=0$, it holds that $x+\mu(1-x)=\mu$; then,
\[
(C|(A \vee B))|((C|A) \wedge (C|B)) = \left\{\begin{array}{ll}
1, &\mbox{if $ABC$ is true,}\\
y+\mu(1-y), &\mbox{if $A\no{B}C$ is true,} \\
\mu, &\mbox{if $\no{AC}$ is true}.
\end{array}
\right.
\]
By coherence, $\mu$ is a linear convex combination of the values $1, y+\mu(1-y)$, that is
\begin{equation}\label{EQ:CONV2}
\mu=\lambda_1+\lambda_2(y+\mu(1-y)) \,,\;\; \lambda_1 \geq 0 \,,\; \lambda_2 \geq 0 \,,\; \lambda_1 + \lambda_2 = 1 \,.
\end{equation}
The equation (\ref{EQ:CONV2}) can be written as 
$\mu(\lambda_1+\lambda_2y) =\lambda_1+\lambda_2y$, 
where $\lambda_1+\lambda_2y>0$; then, $\mu=y+\mu(1-y)=1$ and $(C|(A \vee B))|((C|A) \wedge (C|B)) =1$.
\qed
\end{proof}
\begin{remark}
We observe that 
\[
QC(C|A,C|B)= ((\no{A}\vee C)\wedge(\no{B}\vee C))|(A\vee B)=C|(A\vee B).
\]
Then, the statement of Theorem~\ref{THM:OR} amounts to say that the iterated conditional
$QC(C|A,C|B)|((C|A)\wedge (C|B))$ is equal to 1. This aspect will be analyzed in general in the next section.
\end{remark}
\section{Iterated conditionals and p-entailment}
\label{SECT:MAIN}
In this section we give two results which relate p-entailment and iterated conditioning. In the next result, 
by defining $\F=\{E_1|H_1, E_2|H_2\}$,  $QC(\F)=QC(E_1|H_1, E_2|H_2)$ and $\C(\F)=(E_1|H_1)\wedge(E_2|H_2)$,
we show that, under p-consistency of $\F$,
the iterated conditional
$QC(\F)|(\C(\F))$ is equal to 1.
\begin{theorem}\label{THM:QUASICONG}
	Let  a p-consistent family $\F=\{E_1|H_1, E_2|H_2\}$ be given.
	Then, $QC(\F)|(\C(\F)) = 1$.
\end{theorem}	
	\begin{proof}
	We set 
	\[
	P(E_1|H_1)=x_1 \,,\; P(E_2|H_2)=x_2 \,,\; \prev(\C(\F))=x_{12} \,,\; \prev[\C(\F) \wedge QC(\F)]=\eta \,.	\]
	Moreover, we set $\prev[QC(\F)|\C(\F)]=\mu$.
	Then,
	\[
	QC(\F)|\C(\F) = \C(\F) \wedge QC(\F)+ \mu (1-\C(\F)) \,.  
	\]
	It can be verified that the possible values of the random vector $(\C(\F), \C(\F) \wedge QC(\F))$ are 
	\[
	(1,1) \,,\, (0,0) \,,\, (x_1,x_1) \,,\, (x_2, x_2) \,,\, (x_{12},\eta) \,. 
	\]
	The value $(x_{12},\eta)$ is associated to the constituent $\no{H}_1\no{H}_2$. As we can see, conditionally on $H_1 \vee H_2$ being true, $\C(\F)$ and $\C(\F) \wedge QC(\F)$ coincide; then, by Theorem \ref{EQ-CRQ}, $x_{12}=\eta$, so that $QC(\F) = \C(\F)$. Then,  
	\[
	QC(\F)|\C(\F) = \C(\F) + \mu (1-\C(\F)) = \left\{\begin{array}{ll}
	1, &\mbox{if $\C(\F)=1$,}\\
		\mu, &\mbox{if $\C(\F)=0$,}\\
	x_1+\mu(1-x_1), &\mbox{if $\C(\F) = x_1$,}\\
	x_2+\mu(1-x_2), &\mbox{if $\C(\F) = x_2$,}\\
	x_{12}+\mu(1-x_{12}), &\mbox{if $\C(\F) = x_{12}$.} 
	\end{array}
	\right.
	\]
	By the linearity of prevision, we obtain $\mu=x_{12}+\mu(1-x_{12})$, that is $x_{12}=\mu x_{12}$. Then,  
	\[
	QC(\F)|\C(\F) = \C(\F) + \mu (1-\C(\F)) = \left\{\begin{array}{ll}
	1, &\mbox{if $\C(\F)=1$,}\\
	x_1+\mu(1-x_1), &\mbox{if $\C(\F) = x_1$,}\\
	x_2+\mu(1-x_2), &\mbox{if $\C(\F) = x_2$,}\\
	\mu, &\mbox{if $\C(\F)=0$, or $\C(\F) = x_{12}$.} 
	\end{array}
	\right.
	\]
	We distinguish the following cases: \\ (a) $x_1=x_2=0$; (b) $x_1>0,x_2>0$; (c) $x_1=0,x_2>0$; (d) $x_2=0,x_1>0$. \\
	Case (a). Since $x_1=x_2=0$, it holds that $x_1+\mu(1-x_1)=x_2+\mu(1-x_2)=\mu$, so that  $QC(\F)|\C(\F) \in \{1,\mu\}$. Based on the betting scheme, $\mu=\prev[QC(\F)|\C(\F)]$ is the amount to be paid in order to receive $1$, or $\mu$, according to whether the event $(\C(\F)=1)$ is true, or false, respectively.
	Then, by coherence, it must be the case that $\mu=1$. Therefore, $QC(\F)|\C(\F) = 1$. \\
	Case (b). By coherence, $\mu$ must be a linear convex combination of the values $1$, $x_1+\mu(1-x_1)$, and $x_2+\mu(1-x_2)$, that is,
\begin{equation}\label{EQ:CONV-BIS}
\mu=\lambda_1+\lambda_2(x_1+\mu(1-x_1))+\lambda_3(x_2+\mu(1-x_2)) \,,
\end{equation}
with $\lambda_h \geq 0, h=1,2,3,$ and $\lambda_1+\lambda_2+\lambda_3= 1$.
The equation (\ref{EQ:CONV-BIS}) can be written as
\[
\mu(\lambda_1+\lambda_2x_1+\lambda_3x_2) =\lambda_1+\lambda_2x_1+\lambda_3x_2 \,,
\]
where $\lambda_1+\lambda_2x_1+\lambda_3x_2>0$; then, 
$\mu=x_1+\mu(1-x_1)=x_2+\mu(1-x_2)=1$ and $QC(\F)|\C(\F) = 1$. \\
	Case (c). As $x_1=0$, it holds that $x_1+\mu(1-x_1)=\mu$, so that  
	\[
	QC(\F)|\C(\F) \in \{1,x_2+\mu(1-x_2),\mu\} \,.
	\]
	Then, by coherence, $\mu$ must be a linear convex combination of the values $1,x_2+\mu(1-x_2)$, that is 
	\[
	\mu = \lambda_1 + \lambda_2[x_2+\mu(1-x_2)] \,,\; \lambda_1 + \lambda_2 = 1 \,,\; \lambda_1 \geq 0 \,,\, \lambda_2 \geq 0 \,.
	\]
	It follows that 
	$\mu (\lambda_1 + \lambda_2x_2) = \lambda_1 + \lambda_2x_2$,
	with $\lambda_1 + \lambda_2x_2 > 0$. Then, $\mu=1$  and $QC(\F)|\C(\F) = 1$.  \\
	Case (d). As $x_2=0$, it holds that $x_2+\mu(1-x_2)=\mu$, so that  
	$QC(\F)|\C(\F) \in \{1,x_1+\mu(1-x_1),\mu\}$.
	Then, by coherence, $\mu$ must be a linear convex combination of the values $1,x_1+\mu(1- x_1)$, that is
	\[
	\mu = \lambda_1 + \lambda_2[1_2+\mu(1-x_1)] \,,\; \lambda_1 + \lambda_2 = 1 \,,\; \lambda_1 \geq 0 \,,\, \lambda_2 \geq 0 \,.
	\]
	It follows that
	$\mu (\lambda_1 + \lambda_2x_1) = \lambda_1 + \lambda_2x_1$,
	with $\lambda_1 + \lambda_2x_1 > 0$. Then, $\mu=1$  and $QC(\F)|\C(\F) = 1$. \\
	Therefore, from the p-consistency of the family $\F$ it follows that $QC(\F)|\C(\F) = 1$. 
	\qed
	\end{proof} 

The next theorem shows that the p-entailment of a conditional event $E_3|H_3$ from a p-consistent family $\{E_1|H_1, E_2|H_2\}$ is equivalent to the iterated conditional $(E_3|H_3)|((E_1|H_1)\wedge(E_2|H_2))$ being equal to 1.
\begin{theorem}\label{THM:MAIN}
	Let  three conditional events $E_1|H_1$, $E_2|H_2$, and $E_3|H_3$ be given, where $\{E_1|H_1, E_2|H_2\}$ is p-consistent.
	Then, $\{E_1|H_1, E_2|H_2\}$ p-entails $E_3|H_3$ if and only if
	$(E_3|H_3)|((E_1|H_1)\wedge(E_2|H_2)) = 1$.
\end{theorem}
\begin{proof} 
	$(\Rightarrow)$.
	We observe that by p-consistency 
	$E_1H_1E_2H_2\neq \emptyset$ and then $	(E_1|H_1)\wedge (E_2|H_2)\neq 0$.
	By Theorem \ref{ENTAIL-CS}, $\{E_1|H_1, E_2|H_2\}$ p-entails $E_3|H_3$ if and only if it holds that $QC(\mathcal{S})\subseteq E_3|H_3$ for some $\emptyset \neq  \mathcal{S}\subseteq \{E_1|H_1, E_2|H_2\}$, or $H_3 \subseteq E_3$.
	We observe that, when $H_3 \nsubseteq E_3$, it holds that $\mathcal{S}=\{E_1|H_1\}$, or $\mathcal{S}=\{E_2|H_2\}$, or $\mathcal{S}=\{E_1|H_1,E_2|H_2\}$. 
We  show that  the iterated conditional may be represented as 
	\begin{equation}\label{EQ:REPR-ITER}
	(E_3|H_3)|((E_1|H_1)\wedge(E_2|H_2)) = (E_1|H_1)\wedge(E_2|H_2) + \mu (1 - (E_1|H_1)\wedge(E_2|H_2)) \,,
	\end{equation}		
where $\mu=\prev[(E_3|H_3)|((E_1|H_1)\wedge(E_2|H_2))]$. \\
We distinguish the following  four  cases: \\
	(i) $H_3 \subseteq E_3$; \\
	(ii) $H_3 \nsubseteq E_3$ and $E_1|H_1 \subseteq E_3|H_3$;\\
	(iii) $H_3 \nsubseteq E_3$ and  $E_2|H_2 \subseteq E_3|H_3$; \\
	(iv) $H_3 \nsubseteq E_3$ and  $QC(E_1|H_1,E_2|H_2) \subseteq E_3|H_3$. \\
Case (i). If $H_3 \subseteq E_3$, then $E_3|H_3 = P(E_3|H_3) = 1$. We set $P(E_i|H_i)=x_i,\,  \prev[(E_i|H_i) \wedge (E_j|H_j)]=x_{ij}$ and we recall that 
	\[
	\max \{x_i+x_j-1,0\} \leq x_{ij} \leq \min \{x_i,x_j\} \,.
	\]
	Then, as $x_3=1$, we obtain $x_{13} = x_1,\, x_{23}=x_2$; it follows that for the random vector $((E_1|H_1) \wedge (E_2|H_2),(E_1|H_1) \wedge (E_2|H_2)\wedge(E_3|H_3))$ the possible values are 
	\[
	(1,1) \,,\; (0,0) \,,\; (x_1,x_1) \,,\; (x_2,x_2) \,,\; (x_{12},x_{12}) \,,\; (x_{12}, x_{123}) \,, 
	\]
	where $x_{123}=\prev[(E_1|H_1)\wedge(E_2|H_2)\wedge (E_3|H_3)] = \mu$. As we can see, conditionally on $H_1 \vee H_2 \vee H_3$ being true, $(E_1|H_1) \wedge (E_2|H_2)$ and $(E_1|H_1) \wedge (E_2|H_2)\wedge(E_3|H_3)$ coincide; then, by coherence, $x_{12} = x_{123}$, so that $(E_1|H_1) \wedge (E_2|H_2)\wedge(E_3|H_3)$ and $(E_1|H_1) \wedge (E_2|H_2)$ coincide. 
	Then, (\ref{EQ:REPR-ITER}) is satisfied.
	\\
	Case (ii). As $E_1|H_1 \subseteq E_3|H_3$, by Proposition \ref{PROP:ANDGN} it holds that $E_1|H_1 \wedge E_3|H_3 = E_1|H_1$ and
	$(E_1|H_1) \wedge (E_2|H_2) \wedge (E_3|H_3)=(E_1|H_1)\wedge(E_2|H_2)$. 	Then, (\ref{EQ:REPR-ITER}) is satisfied.\\
	Case (iii). As $E_2|H_2 \subseteq E_3|H_3$, by Proposition \ref{PROP:ANDGN} it holds that  $E_2|H_2 \wedge E_3|H_3 = E_2|H_2$ and 	$(E_1|H_1) \wedge (E_2|H_2) \wedge (E_3|H_3)=(E_1|H_1)\wedge(E_2|H_2)$. 	Then, (\ref{EQ:REPR-ITER}) is satisfied.\\
Case (iv). 
By taking into account that  $QC(E_1|H_1,E_2|H_2) \subseteq E_3|H_3$, 
 the set of possible values of the random vector 	\[
 ((E_1|H_1) \wedge (E_2|H_2)\,,\, QC (E_1|H_1, E_2|H_2)\,,\,(E_1|H_1) \wedge (E_2|H_2) \wedge (E_3|H_3)),
 \] as shown in   Table~\ref{TAB:TABLE1}, is
\[
\{(1,1,1) \,,\, (0,0,0) \,,\, (x_1,1,x_1) \,,\, (x_2,1,x_2) \,,\, (x_{12},\nu_{12},x_{12}) \,,\, (x_{12},\nu_{12},x_{123})\},
\]
where $x_1=P(E_1|H_1), \; x_2=P(E_2|H_2)$, 
$
x_{12}=P[(E_1|H_1)\wedge (E_2|H_2)]$, $\nu_{12}=\prev[QC(E_1|H_1,E_2|H_2)]$,
$x_{123}=P[(E_1|H_1)\wedge (E_2|H_2)\wedge (E_3|H_3)]$.
\begin{table}
		\[
		\small
		\begin{array}{|l|c|c|c|}
 \hline
		  \hspace{1cm}C_h                         &(E_1|H_1)\wedge(E_2|H_2) &QC(E_1|H_1,E_2|H_2) & (E_1|H_1)\wedge(E_2|H_2)\wedge(E_3|H_3) \\ \hline
 E_1H_1       E_2H_2       E_3H_3      &           1 &1 &            1 \\
 E_1H_1       \no{E}_2H_2  E_3H_3      &           0 &0 &            0 \\
 E_1H_1       \no{E}_2H_2  \no{E}_3H_3 &           0 &0 &            0 \\
 E_1H_1       \no{E}_2H_2  \no{H}_3    &           0 & 0&            0 \\
 E_1H_1       \no{H}_2     E_3H_3      &         x_2 & 1&          x_2 \\
 \no{E}_1H_1  E_2H_2       E_3H_3      &           0 & 0&            0 \\
 \no{E}_1H_1  E_2H_2       \no{E}_3H_3 &           0 &0 &            0 \\
 \no{E}_1H_1  E_2H_2       \no{H}_3    &           0 & 0&            0 \\
 \no{E}_1H_1  \no{E}_2H_2  E_3H_3      &           0 & 0&            0 \\
 \no{E}_1H_1  \no{E}_2H_2  \no{E}_3H_3 &           0 &0 &            0 \\
 \no{E}_1H_1  \no{E}_2H_2  \no{H}_3    &           0 & 0&            0 \\
 \no{E}_1H_1  \no{H}_2     E_3H_3      &           0 & 0&            0 \\
 \no{E}_1H_1  \no{H}_2     \no{E}_3H_3 &           0 &0 &            0 \\
 \no{E}_1H_1  \no{H}_2     \no{H}_3    &           0 & 0&            0 \\
 \no{H}_1     E_2H_2       E_3H_3      &         x_1 & 1&          x_1 \\
 \no{H}_1     \no{E}_2H_2  E_3H_3      &           0 &0 &            0 \\
 \no{H}_1     \no{E}_2H_2  \no{E}_3H_3 &           0 & 0&            0 \\
 \no{H}_1     \no{E}_2H_2  \no{H}_3    &           0 & 0&            0 \\
 \no{H}_1     \no{H}_2     E_3H_3      &      x_{12} &\nu_{12} &       x_{12} \\
 \no{H}_1     \no{H}_2     \no{H}_3    &      x_{12} &\nu_{12} &            x_{123} \\
 \hline  
		\end{array}
		\]
\caption{Possible values of the random vector
$((E_1|H_1) \wedge (E_2|H_2)$, $QC (E_1|H_1, E_2|H_2)$, $(E_1|H_1) \wedge (E_2|H_2) \wedge (E_3|H_3))$,
under the assumption that    $QC(E_1|H_1,E_2|H_2) \subseteq E_3|H_3$.
}
		\label{TAB:TABLE1}
\end{table}
As we  can see, conditionally on $H_1 \vee H_2 \vee H_3$ being true (i.e., $\no{H}_1\no{H}_2\no{H}_3$ being false), $(E_1|H_1) \wedge (E_2|H_2)$ and $(E_1|H_1) \wedge (E_2|H_2) \wedge (E_3|H_3)$ coincide; then, by Theorem~\ref{EQ-CRQ} it holds that $x_{12}=x_{123}$, so that  
$(E_1|H_1) \wedge (E_2|H_2) \wedge (E_3|H_3)=(E_1|H_1) \wedge (E_2|H_2)$. Then, (\ref{EQ:REPR-ITER}) is satisfied.\\

Now, by using the representation (\ref{EQ:REPR-ITER}), for the iterated conditional we obtain
\begin{equation}\label{EQ:REPR1}
	(E_3|H_3)|((E_1|H_1)\wedge(E_2|H_2)) = 
	\left\{\begin{array}{ll}
	1, &\mbox{if $E_1H_1E_2H_2$ is true,}\\
	\mu, &\mbox{if $\no{E}_1H_1 \vee \no{E}_2H_2$ is true,}\\
	x_1+\mu(1-x_1), &\mbox{if $\no{H}_1E_2H_2$ is true,}\\
	x_2+\mu(1-x_2), &\mbox{if $E_1H_1\no{H}_2$ is true,}\\
	x_{12}+\mu(1-x_{12}), &\mbox{if $\no{H}_1\no{H}_2$ is true.} 
	\end{array}
	\right.
\end{equation}
Moreover, by  the linearity of prevision it holds that
		\[
		\mu=\prev[(E_3|H_3)|((E_1|H_1)\wedge(E_2|H_2))]
		=x_{12}+\mu(1-x_{12}) \,;
		\]
		from which it follows that $x_{12}=\mu x_{12}$.  
Then, (\ref{EQ:REPR1}) becomes 
\begin{equation}\label{EQ:REPR2}
(E_3|H_3)|((E_1|H_1)\wedge(E_2|H_2))
	= \left\{\begin{array}{ll}
	1, &\mbox{if $E_1H_1E_2H_2$ is true,}\\
	x_1+\mu(1-x_1), &\mbox{if $\no{H}_1E_2H_2$ is true,}\\
	x_2+\mu(1-x_2), &\mbox{if $E_1H_1\no{H}_2$ is true,}\\
	\mu, &\mbox{if $\no{H}_1\no{H}_2 \vee \no{E}_1H_1 \vee \no{E}_2H_2$ is true.} 
	\end{array}
	\right.
\end{equation}
In order to prove that $(E_3|H_3)|((E_1|H_1)\wedge(E_2|H_2)) = 1$, as already done in the proof of Theorem~\ref{THM:OR}, we distinguish the following cases: (a) $x_{12}>0$; (b) $x_{12}=x_1=x_2=0$; (c) $x_{12}=x_1>0,x_2>0$; (d) $x_{12}=x_2=0,x_1>0$; (e) $x_{12}=x_1=0,x_2>0$. \\
	Case (a). As $x_{12}>0$ and $x_{12}=\mu x_{12}$, it follows that $\mu=1$ and then  $x_1+\mu(1-x_1)=x_2+\mu(1-x_2)=1$. Therefore, $(E_3|H_3)|((E_1|H_1)\wedge(E_2|H_2)) = 1$. \\
	Case (b). As $x_1=x_2=0$, it holds that $x_1+\mu(1-x_1)=x_2+\mu(1-x_2)=\mu$, so that  $(E_3|H_3)|((E_1|H_1)\wedge(E_2|H_2)) \in \{1,\mu\}$. We observe that, based on the metaphor of the  betting scheme, $\mu=\prev[(E_3|H_3)|((E_1|H_1)\wedge(E_2|H_2))]$ is the amount to be paid in order to receive $1$, or $\mu$, according to whether $E_1H_1E_2H_2$ is true, or false, respectively.
	Then, by discarding the case where it is received back  what has been paid, coherence requires that $\mu=1$. Therefore $(E_3|H_3)|((E_1|H_1)\wedge(E_2|H_2)) = 1$. \\
	Case (c). By coherence, $\mu$ must be a linear convex combination of the values $1$, $x_1+\mu(1-x_1)$, and $x_2+\mu(1-x_2)$, that is,
\begin{equation}\label{EQ:CONV-TRIS}
\mu=\lambda_1+\lambda_2(x_1+\mu(1-x_1))+\lambda_3(x_2+\mu(1-x_2)) \,,
\end{equation}
with $\lambda_h \geq 0, h=1,2,3,$ and $\lambda_1+\lambda_2+\lambda_3= 1$.
The equation (\ref{EQ:CONV-TRIS}) can be written as
\[
\mu(\lambda_1+\lambda_2x_1+\lambda_3x_2) =\lambda_1+\lambda_2x_1+\lambda_3x_2 \,,
\]
where $\lambda_1+\lambda_2x_1+\lambda_3x_2>0$; then, 
$\mu=x_1+\mu(1-x_1)=x_2+\mu(1-x_2)=1$ and $(E_3|H_3)|((E_1|H_1) \wedge (E_2|H_2)) = 1$. \\
	Case (d). As $x_2=0$, it holds that $x_2+\mu(1-x_2)=\mu$, so that  \[
	(E_3|H_3)|((E_1|H_1)\wedge(E_2|H_2)) \in \{1,x_1+\mu(1-x_1),\mu\} \,.
	\]
	Then, by coherence, $\mu$ must be a linear convex combination of the values $1,x_1+\mu(1-x_1)$, that is 
	\[
	\mu = \lambda_1 + \lambda_2[1_2+\mu(1-x_1)] \,,\; \lambda_1 + \lambda_2 = 1 \,,\; \lambda_1 \geq 0 \,,\, \lambda_2 \geq 0 \,.
	\]
	It follows that 
	$\mu (\lambda_1 + \lambda_2x_1) = \lambda_1 + \lambda_2x_1$,
	with $\lambda_1 + \lambda_2x_1 > 0$. Then, $\mu=1$  and $(E_3|H_3)|((E_1|H_1)\wedge(E_2|H_2)) = 1$.  \\
	Case (e). As $x_1=0$, it holds that $x_1+\mu(1-x_1)=\mu$, so that  \[
	(E_3|H_3)|((E_1|H_1)\wedge(E_2|H_2)) \in \{1,x_2+\mu(1-x_2),\mu\} \,.
	\]
	Then, by coherence, $\mu$ must be a linear convex combination of the values $1,x_2+\mu(1-x_2)$, that is 
	\[
	\mu = \lambda_1 + \lambda_2[x_2+\mu(1-x_2)] \,,\; \lambda_1 + \lambda_2 = 1 \,,\; \lambda_1 \geq 0 \,,\, \lambda_2 \geq 0 \,.
	\]
	It follows that 
	$\mu (\lambda_1 + \lambda_2x_2) = \lambda_1 + \lambda_2x_2$,
	with $\lambda_1 + \lambda_2x_2 > 0$. Then, $\mu=1$  and $(E_3|H_3)|((E_1|H_1)\wedge(E_2|H_2)) = 1$.  \\
	$(\Leftarrow)$.
	Assume that $(E_3|H_3)|((E_1|H_1)\wedge(E_2|H_2))=1$, so that the unique coherent prevision assessment on $(E_3|H_3)|((E_1|H_1)\wedge(E_2|H_2))$ is $\mu=1$.
	From Remark \ref{REM:COMPOUND} it holds that 
	$x_{123}=\mu x_{12}=x_{12}$. Moreover, $x_{123}\leq x_{3}$ \cite[Equation (8)]{GiSa17} and 
	$x_{12}\geq 	\max\{x_1+x_2-1,\, 0\} $ (see Equation (\ref{LOW-UPPER})). Then,	 it holds that
	\[ 
	\max\{x_1+x_2-1,\, 0\} \; \leq \; x_{12} \; = \; x_{123} \; \leq \; x_3 \,,
	\] and, when $x_1=x_2=1$, it follows that 
	$x_{12} = x_{123} = x_3 = 1$. Therefore,  $\{E_1|H_1,E_2|H_2\}$ p-entails $E_3|H_3$. 
	\qed
	\end{proof}
\begin{remark}
We recall that $\{E_1|H_1,E_2|H_2\}$ p-entails 	$QC(E_1|H_1,E_2|H_2)$ (QAND rule, see, e.g., \citeNP{gilio11ecsqaru,GiSa13IJAR}). Then,
Theorem~\ref{THM:QUASICONG} follows  by applying  Theorem~\ref{THM:MAIN} with $E_3|H_3=QC(E_1|H_1,E_2|H_2)$. 
Similar comments can be made for the inference rules examined in Section~\ref{SECT:RULES}.
\end{remark}
In the examples below we show that if $\{E_1|H_1,E_2|H_2\}$ does not p-entail $E_3|H_3$, the iterated conditional 
$(E_3|H_3)|((E_1|H_1)\wedge(E_2|H_2))$ does not coincide with 1.
\begin{example}[Denial of the antecedent]
We consider the rule where the premise set is $\{\no{A},C|A\}$ and the conclusion is $\no{C}$. As is well known, that Denial of the antecedent is neither logically valid in logic nor p-valid in probability logic. Indeed, by defining $P(\no{A})=x,P(C|A)=y,P(\no{C})=z$, it holds that
\[
P(\no{C})=z=1-P(C)=1-[P(C|A)P(A)+P(C|\no{A})P(\no{A})]=1-y(1-x)-P(C|\no{A})x;
\]
Then, when  $x=y=1$, we obtain $z=1-P(C|\no{A})\in[0,1]$; thus,  $\{\no{A},C|A\}$ does not p-entail $\no{C}$.
Then, by Theorem~\ref{THM:MAIN}, the iterated conditional $\no{C}|(\no{A}\wedge (C|A))$ does not coincide with 1. Indeed, by defining $\prev[\no{C}|(\no{A}\wedge (C|A))]=\mu$, it holds that
\[
\no{C}|(\no{A}\wedge (C|A))=
\no{C}\wedge\no{A}\wedge (C|A)+\mu(1-\no{A}\wedge (C|A))=
 \left\{\begin{array}{ll}
	\mu, &\mbox{if $AC$ is true,}\\
	\mu, &\mbox{if $A\no{C}$ is true,}\\
	\mu(1-y), &\mbox{if $\no{A}C$ is true,}\\
	y+\mu(1-y), &\mbox{if $\no{A}\no{C}$ is true.}\\
	\end{array}
	\right.
\] 
If  $y=1$,  we obtain
\[
\no{C}|(\no{A}\wedge (C|A))=
 \left\{\begin{array}{ll}
	\mu, &\mbox{if $AC$ is true,}\\
	\mu, &\mbox{if $A\no{C}$ is true,}\\
	1, &\mbox{if $\no{A}C$ is true,}\\
	0, &\mbox{if $\no{A}\no{C}$ is true,}\\
	\end{array}
	\right.
\]
with $\mu$ being coherent, for every $\mu \in[0,1]$. Therefore, $\no{C}|(\no{A}\wedge (C|A))\neq 1$.
\end{example}
\begin{example}[Affirmation of the consequent]
We consider the rule where the premise set is $\{C,C|A\}$ and the conclusion is $A$. Affirmation of the consequent is neither logically valid in logic nor p-valid in probability logic. Indeed, by defining $P(C)=x,P(C|A)=y,P(A)=z$, and $P(C|\no{A})=t$, it holds that
\[
P(C)=x=P(C|A)P(A)+P(C|\no{A})P(\no{A})=yz+t(1-z).
\]
Then, when  $x=y=1$, we obtain $1=z+t-zt$, that is $z(1-t)=(1-t)$. Therefore, when $t<1$, it follows that  $z=1$. In other words, by adding the premise $P(C|\no{A})<1$  \cite<i.e. what we introduced as a \emph{negated default} in>{gilio16}, it holds that 
\[
P(C)=1, P(C|A)=1, P(C|\no{A})<1 \Rightarrow P(A)=1.
\]
But in general (where no assumptions are made about  $P(C|\no{A})$), $z\in[0,1]$; thus p-entailment of $A$ from  $\{C,C|A\}$ does not hold. Then, by Theorem~\ref{THM:MAIN}, the iterated conditional $A|(C\wedge (C|A))$ does not coincide with 1. 
Indeed, by defining $\prev[A|(C\wedge (C|A))]=\mu$, it holds that
\[
A|(C\wedge (C|A))=
A\wedge C\wedge (C|A)+\mu(1-C\wedge (C|A))=
 \left\{\begin{array}{ll}
	1, &\mbox{if $AC$ is true,}\\
		\mu(1-y), &\mbox{if $\no{A}C$ is true,}\\
	\mu, &\mbox{if $\no{C}$ is true.}\\
	\end{array}
	\right.
\] 
If  $y=1$,  we obtain
\[
A|(C\wedge (C|A))=
 \left\{\begin{array}{ll}
	1, &\mbox{if $AC$ is true,}\\
		0, &\mbox{if $\no{A}C$ is true,}\\
	\mu, &\mbox{if $\no{C}$ is true.}\\
	\end{array}
	\right.
\] 
with $\mu$ being coherent, for every $\mu \in[0,1]$. Therefore, $A|(C\wedge (C|A))\neq 1$.
\end{example}
As another example, we could consider  \emph{Transitivity}, where $\{C|B,B|A\}$ is the premise set and $C|A$ is the conclusion. The p-entailment does not hold, indeed the assessment $(1,1,z)$ on $\{C|B,B|A,C|A\}$ is coherent for any $z\in[0,1]$. Then, by Theorem \ref{THM:MAIN}, the iterated conditional $(C|A)|(C|B)\wedge (B|A)$ does not coincide with 1. But, by adding  the negated default  $P(\no{A}|(A\vee B))<1$  it holds that \cite[Theorem 5]{gilio16}
\[
P(C|B)=1, P(B|A)=1, P(\no{A}|(A\vee B))<1\, \Rightarrow\, P(C|A)=1.
\]
\section{Concluding remarks}
The results of this paper are based on the notions of conjoined conditionals  and iterated  conditionals. These  objects, introduced in recent papers by Gilio and Sanfilippo, are defined  in the setting of coherence by means of suitable conditional random quantities with values in the interval $[0,1]$. 
By exploiting the logical implication of Goodman and Nguyen, we have shown that  $A|H$ p-entails $B|K$ if and only if $(B|K)|(A|H) = 1$. Moreover, we have shown that a p-consistent family  $\F=\{E_1|H_1,E_2|H_2\}$ p-entails a conditional event $E_3|H_3$ if and only if  $E_3|H_3=1$, or  $(E_3|H_3)|QC(\S)=1$ for some nonempty subset $\S$ of $\F$. We have also applied our result considered the inference rules And, Cut, Cautious Monotonicity, and  Or of System~P  and the inference rules  Modus Ponens, Modus Tollens, and  Bayes.  We have also shown that the iterated conditional $QC(\F)|\C(\F)$ is equal to 1 for every p-consistent  family $\F=\{E_1|H_1,E_2|H_2\}$.
Then, we have characterized the p-entailment of  $E_3|H_3$ from a p-consistent family $\F$ by showing that it amounts to
the condition $(E_3|H_3)|\C(\F)=1$.
Finally, we examined two examples (Denial of the Antecedent and Affirmation of the Consequent)  when the  p-entailment of the conditional event $E_3|H_3$ from a p-consistent family $\{E_1|H_1,E_2|H_2\}$ does not hold by also showing that $(E_3|H_3)|((E_1|H_1)\wedge(E_2|H_2)   )\neq 1$. Concerning  the  Affirmation of the Consequent, we also showed that (a kind of conditional) p-entailment holds if we add
a suitable negated default in the set of premises. Psychologically, this could serve as a new explanation why some people interpret Affirmation of the Consequent as a valid argument form. Indeed, this argument form plays an important r\^{o}le in abductive reasoning in philosophy of science (e.g., where conclusions about possible causes/diseases are derived from effects/symptoms).  Future work is needed to  explore such applications of the presented theory and to explore further formal desiderata  also related to the deduction theorem.

%
\end{document}